\documentclass[10pt]{article}
\usepackage[T1]{fontenc}
\usepackage{amsmath,amsfonts,amsthm,mathrsfs,amssymb}
\usepackage{cite}
\usepackage{graphicx,float}
\usepackage{subfigure}
\usepackage{placeins}
\usepackage{color}
\usepackage{indentfirst}
\usepackage[bookmarks=true]{hyperref}
\numberwithin{equation}{section}
\newcommand{\R}{{\mathbb R}}

\newcommand{\be}{\begin{eqnarray}}
\newcommand{\ben}{\begin{eqnarray*}}
\newcommand{\en}{\end{eqnarray}}
\newcommand{\enn}{\end{eqnarray*}}

\newcommand{\pa}{\partial}

\newcommand{\curl}{{\rm curl\,}}

\newcommand{\divv}{{\rm div\,}}

\newcommand{\G}{\Gamma}

\definecolor{rot}{rgb}{0.000,0.000,0.000}
\definecolor{rot1}{rgb}{0.000,0.000,0.000}
\begin{document}
\renewcommand{\theequation}{\arabic{section}.\arabic{equation}}
\begin{titlepage}
\title{Weighted integral solvers for elastic scattering by open arcs in two dimensions}

\author{Oscar P. Bruno\thanks{Department of Computing \& Mathematical Sciences, California Institute of Technology, 1200 East California Blvd., CA 91125, United States. Email:{\tt obruno@caltech.edu}}\;,
Liwei Xu\thanks{School of Mathematical Sciences, University of Electronic Science and Technology of China, Chengdu, Sichuan 611731, China. Email: {\tt xul@uestc.edu.cn}}\;,
Tao Yin\thanks{Department of Computing \& Mathematical Sciences, California Institute of Technology, 1200 East California Blvd., CA 91125, United States. Email:{\tt taoyin89@caltech.edu}}}
\end{titlepage}
\maketitle

\begin{abstract}
  We present a novel approach for the numerical solution of problems
  of elastic scattering by open arcs in two dimensions. Our
  methodology relies on the composition of weighted versions of the
  classical operators associated with Dirichlet and Neumann boundary
  conditions in conjunction with a certain "open-arc elastic
  Calder\'on relation" whose validity is demonstrated in this paper on
  the basis of numerical experiments, but whose rigorous mathematical
  proof is left for future work. Using this Calder\'on relation in
  conjunction with spectrally accurate quadrature rules and the
  Krylov-subspace linear algebra solver GMRES, the proposed overall
  open-arc elastic solver produces results of high accuracy in small
  number of iterations---for low and high frequencies alike. A variety
  of numerical examples in this paper demonstrate the accuracy and
  efficiency of the proposed methodology.  {\bf Keywords:} Elasticity,
  open arc, Calder\'on relation, second-kind integral solver.
\end{abstract}

\section{Introduction}
\label{sec:1}

We consider the problem of numerical evaluation of elastic waves
diffracted by {\em infinitely thin open
  surfaces~}\cite{MG14,PG15,PG18,PGH17,W06}.  This problem plays
fundamental roles in a number of important applications in science and
engineering, including non-destructive testing of materials,
characterization of fractures, energy production from natural gas and
geothermal resources, mining, etc. These problems present considerable
mathematical and computational challenges in view of the highly
oscillatory character of the associated time-harmonic elastic fields
and, for solvers based on use of volumetric discretizations, the
unbounded character of the physical domains that must be considered in
connection with the aforementioned applications.

As is the case for wave scattering problems in acoustics and
electromagnetics~\cite{CK98,N01}, the boundary integral equation
methods in elasticity require discretization of domains of lower
dimensionality than those required by volumetric discretization
methods (such as finite difference or finite element
methods~\cite{BHSY,PGH17}).  These equations can generally be treated
effectively even for high frequencies~\cite{BK01,BG18,BLR14,BL12}. As
is well known, however, the classical boundary integral equations for
open surfaces (or, in two dimensions, open arcs) are not second-kind
Fredholm integral equations. This setting presents some
difficulties. On one hand, as the eigenvalues of the left-hand side
operators accumulate around zero and/or infinity, solution of these
problems by means of Krylov-subspace iterative solvers such as GMRES
(which are commonly-used in conjucntion with accelerated
integral-equation solvers) often requires large numbers of iterations
for convergence, and thus, large computing costs. Additionally, the
evaluation of the Hadamard finite part of the hypersingular integral
operator associated to the elastic Neumann problem~\cite{AD95} has
also remained a significant challenges in this context~\cite{BLR14}.

These difficulties were recently addressed, for the 2D and 3D acoustic
and 2D electromagnetic contexts, in the
contributions~\cite{BL12,BL13,LB15}---which, in particular, introduced
Fredholm integral equations of the second-kind, and associated
numerical algorithms, for problems of acoustic scattering by open arcs
and surfaces. The second-kind equations were obtained in these
articles by utilizing compositions of appropriately modified versions
$\textcolor{rot}{S_a^w}$ and $\textcolor{rot}{N_a^w}$ (which incorporate explicitly the singular character of the integral equation densities) of the acoustic single-layer and
hypersingular integral operators $\textcolor{rot}{S_a}$ and $\textcolor{rot}{N_a}$. As shown in~\cite{BL12,LB15}, a generalization of the classical closed-surface
Calder\'on formulas holds for the $\textcolor{rot}{S_a^w}$ and $\textcolor{rot}{N_a^w}$ operators in the open arc case, which, in particular, gives rise to the aforementioned
second-kind integral formulations for the acoustic open-arc
problems. It was verified numerically in these contributions that, as
predicted by theory, the eigenvalues of the proposed second-kind
operators remain bounded away from zero and infinity, even for
problems of high frequency.

Multiple challenges arise as extensions of these methods to elastic
open-arc problems are attempted. At a basic level, elastic Calder\'on
formulas have not been studied even in the closed-surface
case---likely, on account of the fact that, unlike the acoustic wave
case, the classical Neumann-Poincar\'e double-layer operator $K$ and
its adjoint $K^*$, which play important roles in the Calder\'on
relations, are not compact in the elastic
case~\cite{AJKKY,KGBB79}. (Reference~\cite{DSS15} mentions all the
elastic operators relevant to the closed-surface Calder\'on calculus
but it does not utilize the Calder\'on projectors as regularization
tools). But recent results~\cite{AJKKY} have established that the
closed surface \textcolor{rot}{elastostatic} ($\omega=0$) double-layer
integral operator is polynomially compact, which suggests that the
composition $NS$ of closed-surface single-layer and hypersingular
integral operators may be a Fredholm operator of second-kind even in
the case $\omega\ne 0$.  This in fact established in
Section~\ref{sec:3.2} below, where it is additionally shown that the
eigenvalues of the operator $NS$ accumulate at a certain point which
depends on the elastic Lam\'e parameters. Moreover, for the Dirichlet
problem particulary, an artificial traction operator $\widetilde{T}$
can be introduced (see (\ref{astress}) for which the resulting
double-layer operator $\widetilde{K}$ is
compact~\cite{H98,KGBB79}---which results in an alternative
second-kind Calder\'on-like formula for the composition
$\widetilde{N}S$. But this formula is only applicable for Dirichlet
problem since this artificial traction operator does not correspond to
an actual physical traction.

Relying on the newly studied property of the closed-surface Calder\'on
formula for elastic wave, it is naturally to consider the extension of
the existing second-kind Fredholm integral equations for acoustic
open-arcs~\cite{BL12,LB15} to the elastic case. In view of these
contributions, and unlike the
approaches~\cite{AS91,CN00,JR04,M96,PS60,YS88}, we consider the
composition of weighted versions of the classical single-layer and
hypersingular integral operators~\cite{AD95}. We find that the
benefits of this approach are two-fold. On one hand, the new method
enjoys high-order accuracy: the weighted versions $S^w$, $N^w$ of the
single-layer and hypersingular operators extract the solutions' edge
singularity~\cite{CDD03} explicitly and the applied quadrature rules
provide spectral convergence. On the other hand, the method gives rise
to well-behaved linear algebra: as numerically demonstrated in Section
\ref{sec:3.3}, \textcolor{rot}{the eigenvalues of $N^wS^w$ and
  $\widetilde{N}^wS^w$ are bounded away from zero and infinity. (A
  theoretical proof of this fact is left for future work.)}. And, as
desired, the composite operator $N^wS^w$ or $\widetilde{N}^wS^w$
requires small number of iterations when used in conjunction with the
linear iterative solver GMRES. The new formulation for the Neumann
problem is especially beneficial, as it give rise to
order-of-magnitude improvements in computing times over the
corresponding hypersingular formulation. Such gains do not occur for
the Dirichlet problem, although the new formulation
$\widetilde{N}^wS^w$ requires fewer iterations than $S^w$, since the
application of the operator $S^w$ is significantly less expensive than
the application of the operator $\widetilde{N}^w$.

A number of additional techniques are proposed to accurately evaluate
the elastic hypersingular integral operators $N^w$ and
$\widetilde{N}^w$. For closed surface scattering problems in
elasticity, the novel and exact regularized formulation presented in
\cite{BXY17,YHX17} show that the hypersingular operator in two
dimensions can be transformed into a composition of weakly-singular
integrals and tangential-derivative operators that involve the
G\"unter derivative and integration-by-parts. For the weighted
hypersingular operators $N^w$ and $\widetilde{N}^w$ in the open-arc
case, thanks to the edge-vanishing weight function $w$, the results in
the closed-surface case can be extended to the open arc case since all
singular terms arising from the integration-by-parts calculation can
be eliminated. Additionally, the tangential derivative evaluations,
which are approximated in our method by means of FFTs, has proven more
efficient than the alternative treatment~\cite{K95}.

The remainder of this paper is organized as follows. Section
\ref{sec:2} describes the Dirichlet and Neumann problems of elastic
scattering by open surfaces, along with the associated classical
boundary integral formulations and including a brief discussion of
certain associated challenges. Section \ref{sec:3.1} introduces new
weighted operators that explicitly account for the known singular
character of the solutions, and Section \ref{sec:3.2} presents our
novel investigation of the closed-surface Calder\'on formula,
including numerical verifications of a polynomial-compactness
theoretical result. Section~\ref{sec:3.3} then presents a numerical
examination of the Calder\'on relation for open arcs in light of the
eigenvalue distributions obtained for non-trivial open arc problems.
An exact regularized formulation for the hyper-singular operator is
presented in Section \ref{sec:3.4}. The high order quadrature rules we
use for evaluation of the new integral operators are described in
Section~\ref{sec:4.1}. Numerical demonstrations presented in
Section~\ref{sec:4.2}, for both low and high frequencies and for
various geometries, demonstrate the high-accuracy and high-order of
convergence enjoyed by the proposed approach, as well as the reduced
numbers of GMRES linear-algebra iterations required by the algorithm
for convergence.

\section{Preliminaries}
\label{sec:2}
\subsection{The elastic scattering problem}
\label{sec:2.1}

Let $\Gamma$ denote a smooth open arc in the plane $\mathbb{R}^2$. The
complement of $\Gamma$ in $\mathbb{R}^2$ is occupied by a linear
isotropic and homogeneous elastic medium characterized by the Lam\'e
constants $\lambda,\mu$ satisfying $\mu>0$, $\lambda+\mu>0$, and by
the mass density $\rho>0$. Suppressing the time-harmonic dependence
$e^{-i\omega t}$ in which $\omega>0$ is the frequency, the
displacement field is the solution of the time-harmonic Navier
equation \be
\label{navier}
\Delta^*u+\rho\omega^2u=0 \quad\mbox{in}\quad\R^2\backslash\Gamma, \en
together the appropriate Kupradze radiation condition~\cite{KGBB79} at
infinity. On $\Gamma$ the solution is assumed to satisfy either the
Dirichlet boundary condition \ben u=F\quad\mbox{on}\quad\Gamma  \enn
or the Neumann boundary condition \ben \quad T(\pa,\nu)u=G
\quad\mbox{on}\quad\Gamma.  \enn Here $\Delta^{*}$ is the Lam\'e
operator defined by \ben
\label{LameOper}
\Delta^* = \mu\,\mbox{div}\,\mbox{grad} + (\lambda + \mu)\,\mbox{grad}\, \mbox{div}\,,
\enn
and $T=T(\pa,\nu)$ is the traction (or stress) operator on the boundary defined as
\be
\label{stress}
Tu:=2 \mu \, \partial_{\nu} u + \lambda \,
\nu \, \divv u+\mu \nu\times \curl u,\quad \nu=(\nu^1,\nu^2){^\top},
\en
in which $\nu$ is the unit normal to the boundary $\G$ and $\partial_\nu:=\nu\cdot\nabla$ is the normal derivative.

Along with~(\ref{stress}), we define the modified (unphysical)
traction operator \be
\label{astress}
\widetilde{T}u:=(\mu+\widetilde{\mu}) \, \partial_{\nu} u +
\widetilde{\lambda} \, \nu \, \divv u+\widetilde{\mu}\nu\times \curl
u, \en where $\widetilde{\lambda}+\widetilde{\mu}=\lambda+\mu$. It
easily follows that $\widetilde{T}=T$ when
$\widetilde{\lambda}=\lambda, \widetilde{\mu}=\mu$. The selection
\be\label{tildes} \widetilde{\mu}=\frac{\mu(\lambda+\mu)}{\lambda+3\mu},
\quad \widetilde{\lambda}=\lambda+\mu-\widetilde{\mu}, \en which is
made throughout this article, is justified in~\ref{sec:3.2} and
subsequent sections.

\subsection{Boundary integral equations}
\label{sec:2.2}

It follows from potential theory that the solution of (\ref{navier})
under Dirichlet and Neumann boundary conditions can be expressed in
terms of single- and double-layer potentials,
 \be
\label{DirichletS}
u(x)  = (\mathcal{S}\Psi)(x):=  \int_{\Gamma}E(x,y)) \phi(y)\,ds_y, \quad \forall\,x\in\R^2\backslash\Gamma,
\en
and
\be
\label{NeumannD}
u(x) = (\mathcal{D}\Psi)(x):= \int_{\Gamma}(T(\pa_y,\nu_y)
E(x,y))^\top \psi(y)\,ds_y, \quad \forall\,x\in\R^2\backslash\Gamma,
\en respectively. Here, calling $\gamma_{k_t}(x,y)$  the
fundamental solution of the Helmholtz equation in $\R^2$ with wave
number $k_t$, \be
\label{HelmholtzFS}
\gamma_{k_t}(x,y) = \frac{i}{4}H_0^{(1)}(k_t|x-y|), \quad x\ne y, \en
and letting \ben k_s := \omega/c_p,\quad k_p = \omega/c_s \enn denote
the wave number of the compressional and shear waves for isotropic
dymanic elasticity, respectively, where \ben c_p=\sqrt{\mu/\rho}\quad
\mbox{and}\quad c_s=\sqrt{(\lambda+2\mu)/\rho}, \enn $E(x,y)$ denotes
the fundamental displacement tensor for the Navier equation in $\R^2$:
\ben E(x,y)=\frac{1}{\mu}\gamma_{k_s}(x,y)I+\frac{1}{\rho\omega^2}
\nabla_x\nabla_x^\top
\left[\gamma_{k_s}(x,y)-\gamma_{k_p}(x,y)\right].  \enn As is known,
using the single-layer and hypersingular operators \be
\label{SBIO}
S[\phi](x)  &=& \int_\G E(x,y)\phi(y)ds_y,\quad x\in\G
\en
and
\be
\label{HBIO}
N[\psi](x) &=& T(\pa_x,\nu_x)\int_\G
(T(\pa_y,\nu_y)E(x,y))^\top\varphi(y)ds_y,\quad x\in\G,  \en the
Dirichlet and Neumann problems reduce to the boundary integral
equations \be
\label{BIEe}
S[\phi]=F,\quad N[\psi]=G \quad\mbox{on}\quad\G.
\en

\begin{figure}[ht]
\centering
\begin{tabular}{cc}
\includegraphics[scale=0.20]{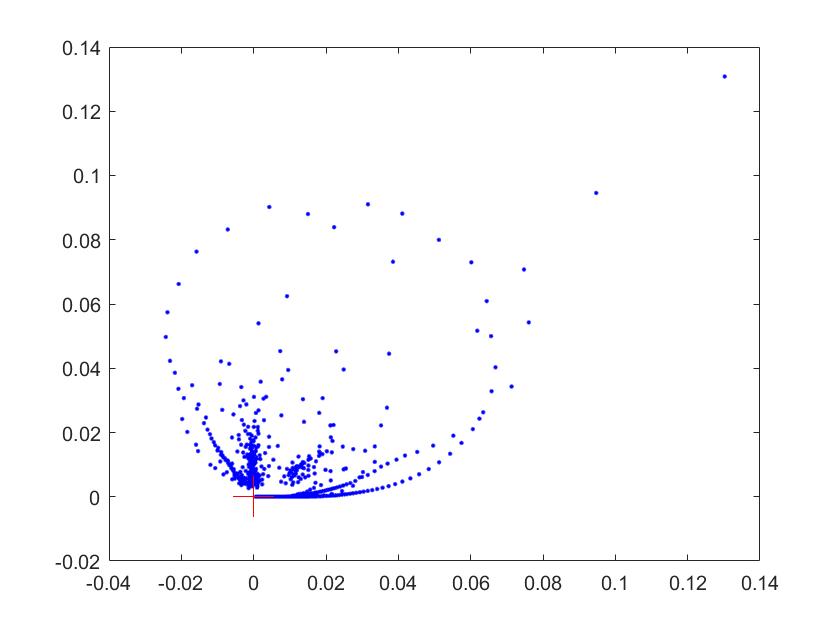} &
\includegraphics[scale=0.20]{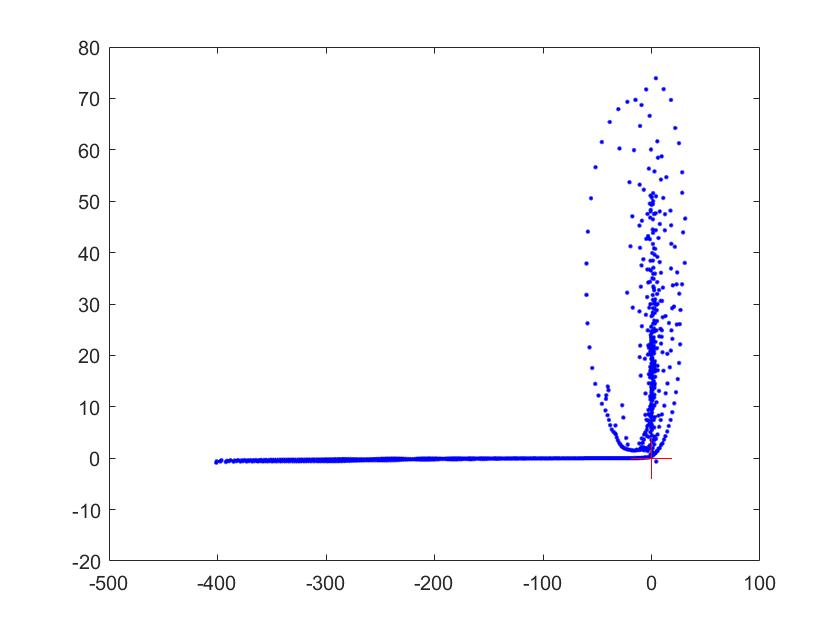} \\
(a) $S^w$ & (b) $N^w$
\end{tabular}
\caption{Eigenvalue distribution of $S^w$ and $N^w$ for the spiral-shaped arc. Red cross: $(0,0)$.}
\label{EeigNS}
\end{figure}

It is well known that, as demonstrated in Fig.\ref{EeigNS}, the
eigenvalues of the single-layer and hypersingular integral operators
in equations (\ref{BIEe}) accumulate at zero and infinity,
respectively. As a result (and as illustrated in numerical tests in
section \ref{sec:4.2}), the solution of these equations by means of
Krylov-subspace iterative solvers such as GMRES generally requires
large numbers of iterations. In addition, as discussed in section
\ref{sec:3.1}, the solutions $\phi$ and $\psi$ of equations
(\ref{BIEe}) are not smooth at the end-points of $\Gamma$, and, thus,
they give rise to low order convergence (and require high
discretization of the densities for a given accuracy) unless such
singularities are appropriately treated.

\section{Weighted operators and new integral solvers}
\label{sec:3}
\subsection{Regularity and singular behavior at the edge}
\label{sec:3.1}

The singular character of the solutions of scattering by open arcs is
well documented~\cite{CDD03}. In particular, it is known that
$\phi,\varphi$ can be expressed in the forms \be
\label{rs1}
\phi\sim \frac{\xi_1}{\sqrt{d}}+\eta_1,\quad \varphi\sim\xi_2\sqrt{d}+\eta_2,
\en
where $d$ denotes the distance to the edge, $\xi_1$ and $\xi_2$ are smooth cut-off functions, and where the functions $\eta_1$ and $\eta_2$ are somewhat smoother than $\phi$ and $\varphi$. If the curve itself and
the boundary functions $F$ and $G$ are infinitely differentiable, it follows that
\be
\label{rs2}
\phi= \frac{\alpha}{\sqrt{d}},\quad \varphi=\beta\sqrt{d},
\en
where $\alpha$ and $\beta$ are infinitely differentiable functions throughout $\Gamma$, up to and including the endpoints. Thus, the singular character of these solutions is fully characterized by the factors $d^{1/2}$ and $d^{-1/2}$ in equation (\ref{rs2}).

In view of the regularity results (\ref{rs2}), we introduce a positive
integral weight $w(x)>0$ with asymptotic behavior $w\sim\sqrt{d}$ at
the edge (by which it is implied that the quotient $w/\sqrt{d}$ is
infinitely differentiable up to the edge), and we define the weighted
operators \ben S^w[\alpha]=S\left[\frac{\alpha}{w}\right],\quad
N^w[\beta]=N[\beta w], \enn so that equations~\eqref{BIEe} for the
Dirichlet and Neumann problems may be expressed in the forms \be
\label{BIE1}
S^w[\alpha]=F,\quad N^w[\beta]=G \quad\mbox{on}\quad\Gamma, \en
respectively. In what follows a smooth parameterizations
$x(t)=(x_1(t),x_2(t))$ defined in the interval $[-1,1]$ is used for
the curve $\Gamma$, and a canonical choice is made for the weight $w$:
$w(x(t))=\sqrt{1-t^2}$.

\subsection{Calder\'{o}n relation for closed surfaces}
\label{sec:3.2}
It is well known that for acoustic scattering by closed smooth
surfaces there holds the Calder\'{o}n formula \ben
\textcolor{rot}{N_{a,c}S_{a,c}=-\frac{I}{4}+(K_{a,c}^*)^2}, \enn where
$\textcolor{rot}{S_{a,c},K_{a,c}^*}$ and $\textcolor{rot}{N_{a,c}}$
are the corresponding single-layer, adjoint of double-layer and
hyper-singular boundary integral operators of acoustic scattering by
closed surfaces, respectively, and $\textcolor{rot}{K_{a,c}^*}$ is
compact. For  elastic scattering by closed surfaces, we can obtain a
similar Calder\'{o}n formula for the single-layer and hyper-singular
boundary integral operators $S,N$. \textcolor{rot}{However, as discussed below, the adjoint  $K^*$
  of the double-layer operator} is non-compact---a fact that
makes the elastic scattering problem more challenging than its acoustic
counterpart. But as indicated in what follows, $NS$ can still be
viewed as a compact perturbation of a multiple of the identity
operator.

The operator $\textcolor{rot}{K^*}$ is given by \be
\label{DBIO}
\textcolor{rot}{K^*}[\psi](x) &=&\int_\G
T(\pa_x,\nu_x)E(x,y)\varphi(y)ds_y,\quad x\in\G.  \en Utilizing the
artificial traction operator (\ref{astress}), we additionally define
the operators \ben
\textcolor{rot}{\widetilde{K}^*}[\psi](x)  &=&\int_\G \widetilde{T}(\pa_x,\nu_x)E(x,y)\varphi(y)ds_y,\quad x\in\G,\\
\widetilde{N}[\psi](x) &=& \widetilde{T}(\pa_x,\nu_x)\int_\G
(\widetilde{T}(\pa_y,\nu_y)E(x,y))^\top\varphi(y)ds_y,\quad x\in\G.
\enn For the special choice~\eqref{tildes} of the constants
$\widetilde{\lambda}$ and $\widetilde{\mu}$, it is
known~\cite{H98,KGBB79} that the kernel of $\widetilde{K}^*$ is
weakly-singular, which implies that that $\widetilde{K}^*$ is a
compact operator. However, the kernel of $K^*$ is strongly singular
and, consequently, the operator $K^*$ itself is not compact. As indicated above, this
presents a singificant difficulty in our context which can, however,
be resolved by appealing to a recent result~\cite{AJKKY}in the
elasto-static context---as indicated in what follows.

\textcolor{rot}{Let $K_0^*$ denote the adjoint of the elastic
  double-layer operator in the zero-frequency case $\omega=0$. The
  spectral properties of this kind of operator are best known in the
  electrostatic case (which involves the Laplace equation), where they
  relate to plasmonics and cloaking by anomalous localized resonances,
  which occur at eigenvalues and at the accumulation point of
  eigenvalues, respectively. The recent contribution~\cite{AJKKY}
  shows that cloaking by anomalous localized resonance also occurs in
  the elastic case. A fundamental result in that paper states that
  $(K_0^*)^2-C_{\lambda,\mu}^2I$ is a compact operator (and, in
  particular, that $K_0^*$ is polynomially compact)}, where
$C_{\lambda,\mu}$ is a constant that depends on the Lam\'e parameters:
\ben C_{\lambda,\mu}=\frac{\mu}{2(\lambda+2\mu)}.  \enn

\textcolor{rot}{This result is useful in our context. Indeed, noting
  that $K^*-K_0^*$ has a weakly-singular kernel it follows that
  $K^*-K_0^*$ is a compact operator, and we obtain \ben
  (K^*)^2-C_{\lambda,\mu}^2I= (K^*-K_0^*)(K^*+K_0^*)
  +(K_0^*)^2-C_{\lambda,\mu}^2I \enn is compact.} We thus obtain the
elastic Calder\'{o}n formulae for closed smooth surfaces \be
\label{calderonclosed}
NS=-\left(\frac{1}{4}-C_{\lambda,\mu}^2\right)I+K^{(1)}, \quad
\widetilde{N}S=-\frac{1}{4}I+K^{(2)}, \en where $K^{(1)},K^{(2)}$ are
compact operators.

\begin{figure}[htbp]
\centering
\begin{tabular}{cc}
\includegraphics[scale=0.2]{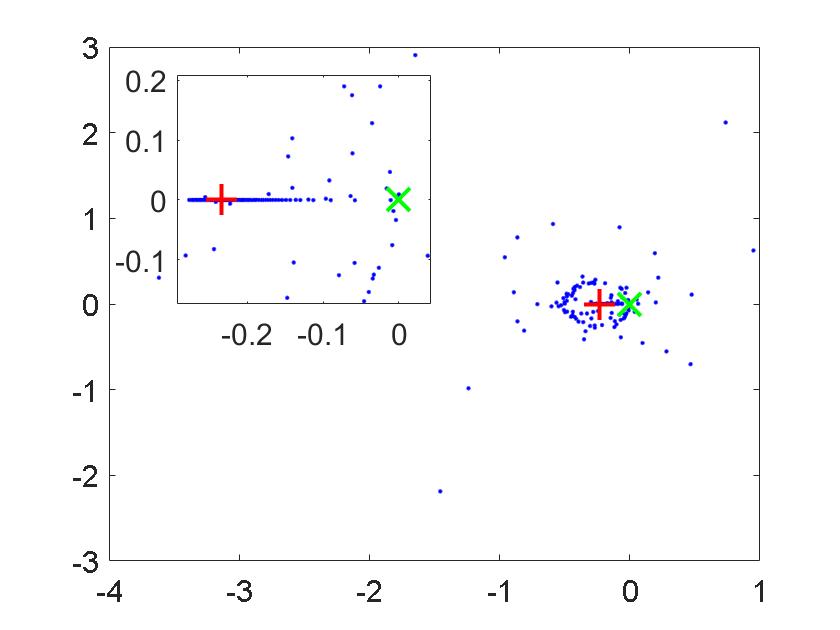} &
\includegraphics[scale=0.2]{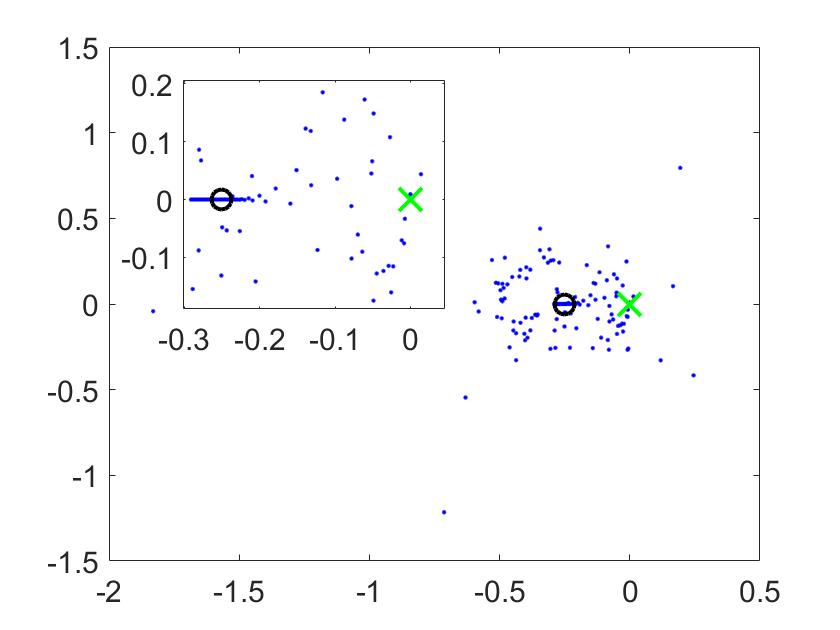} \\
(a) $NS$, $\lambda=2$ & (b) $\widetilde{N}S$, $\lambda=2$\\
\includegraphics[scale=0.2]{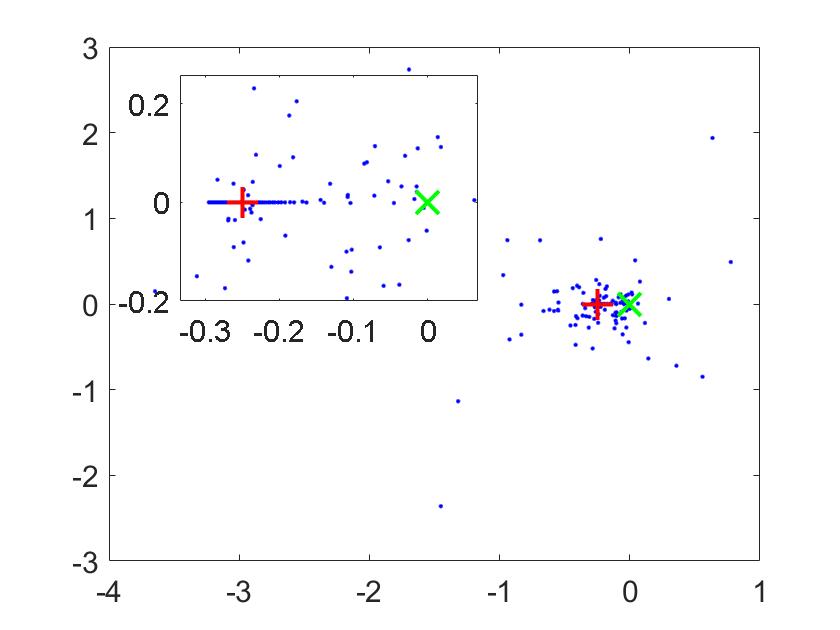} &
\includegraphics[scale=0.2]{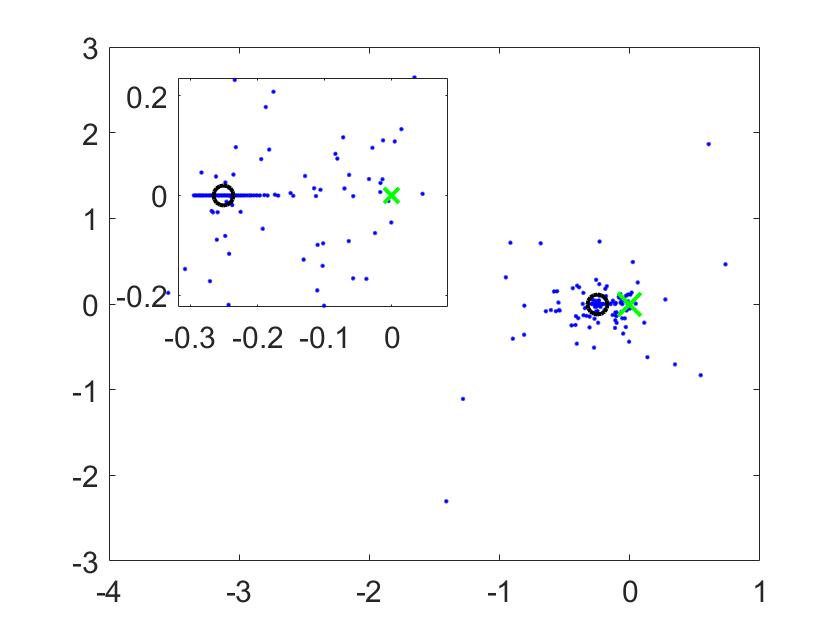} \\
(c) $NS$, $\lambda=100$ & (d) $\widetilde{N}S$, $\lambda=100$\\
\includegraphics[scale=0.2]{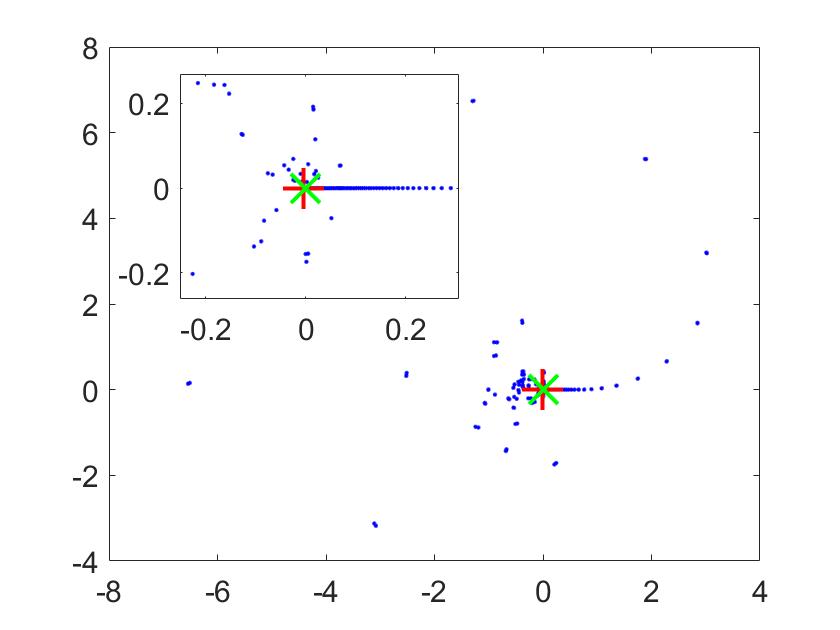} &
\includegraphics[scale=0.2]{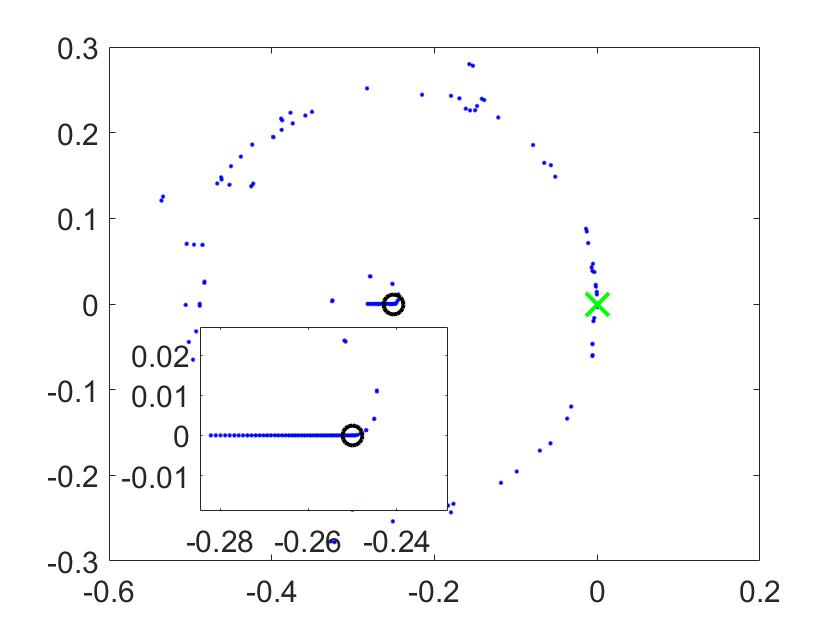} \\
(e) $NS$, $\lambda=-0.99$ & (f) $\widetilde{N}S$, $\lambda=-0.99$
\end{tabular}
\caption{Eigenvalue distribution of $NS$ and $\widetilde{N}S$ for a
  circular scatterer. Red points (+): $(-1/4+C_{\lambda,\mu}^2, 0)$; Black
  points ($\circ$): $(-1/4,0)$; Green points ($\times$): $(0,0)$. As predicted by theory, the
  eigenvalues of $NS$ and $\widetilde{N}S$ accumulate at the red and
  black points, respectively, and they are bounded away from zero and infinity.}
\label{CloseEeigNS}
\end{figure}

To verify the above results numerically, we consider the problem of
elastic scattering by a circular scatterer of radius one, and we
choose $\rho=1$, $\mu=1$, $\omega=50$.  Fig.~\ref{CloseEeigNS}
displays the eigenvalue distribution of $NS$ and $\widetilde{N}S$ for
various values of $\lambda$. (The eigenvalue computation was based on
the exact regularized formulations for the hypersingular integral
operator given in section \ref{sec:3.4}, implemented by means of the
high-order Nystr\"om methodology~\cite{CK98} together with FFT for
evaluation of tangential derivatives.) Choosing a large enough number
of descritization points the eigenvalues of $NS$ and $\widetilde{N}S$
are seen to accumulate at $(-1/4+C_{\lambda,\mu}^2, 0)$ and
$(-1/4,0)$, respectively, as predicted by the result embodied in
Equation~\eqref{calderonclosed}.

\begin{figure}[htbp]
\centering
\begin{tabular}{cc}
\includegraphics[scale=0.2]{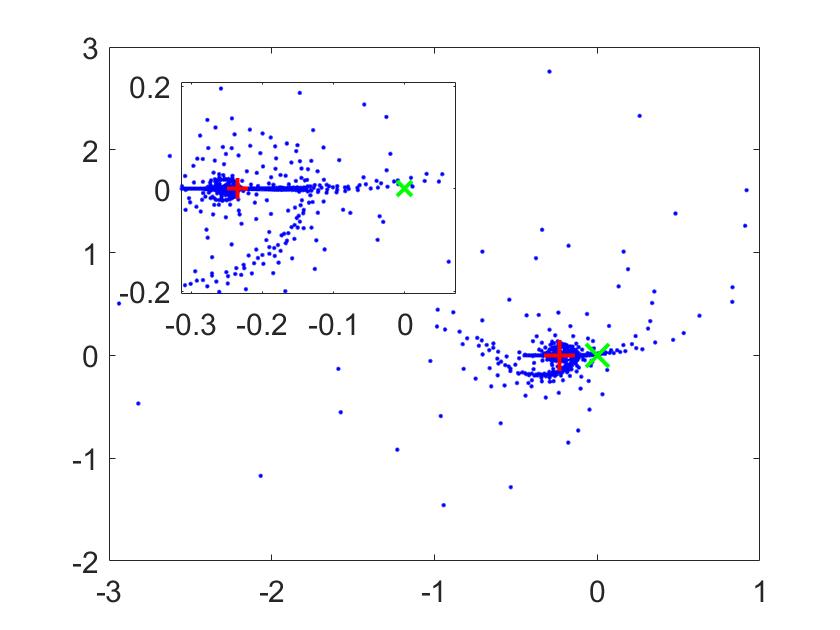} &
\includegraphics[scale=0.2]{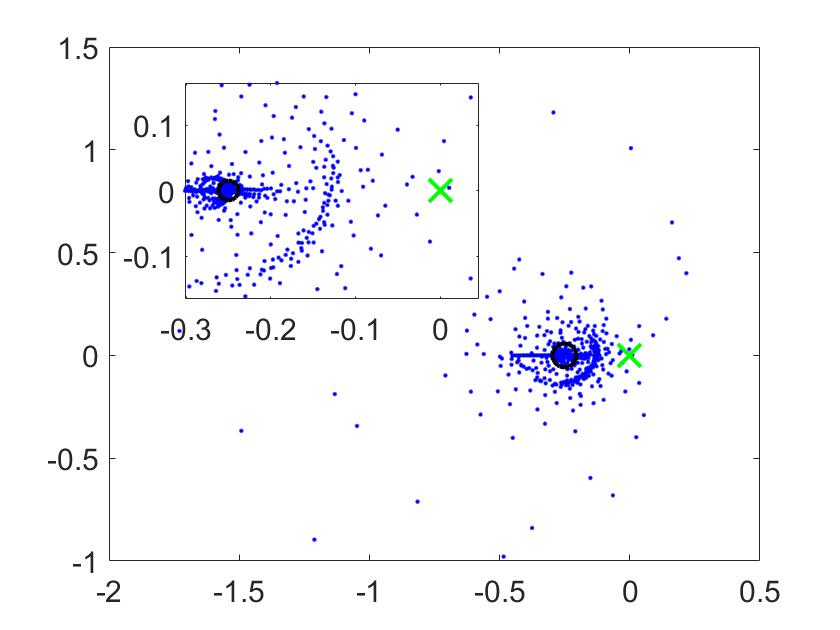} \\
(a) $N^wS^w$, $\lambda=2$ & (b) $\widetilde{N}^wS^w$, $\lambda=2$\\
\includegraphics[scale=0.2]{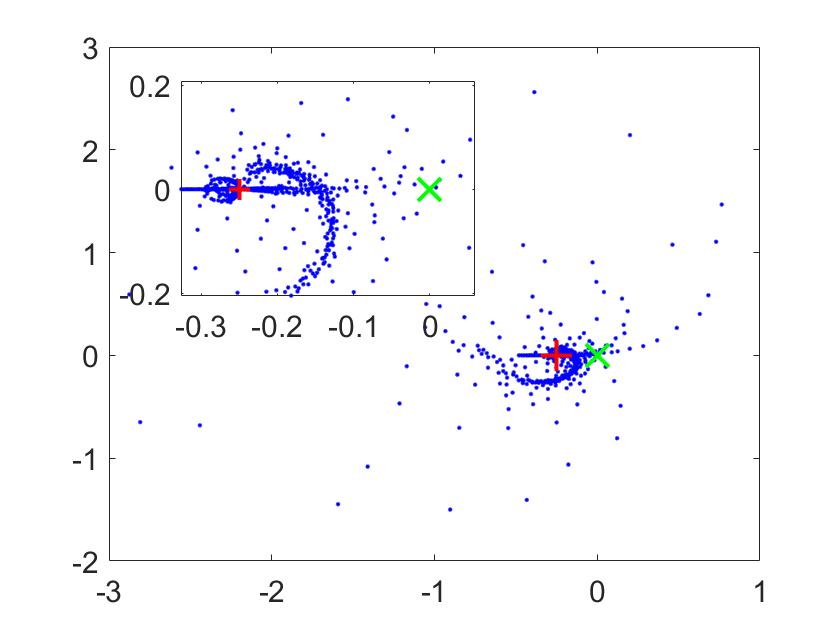} &
\includegraphics[scale=0.2]{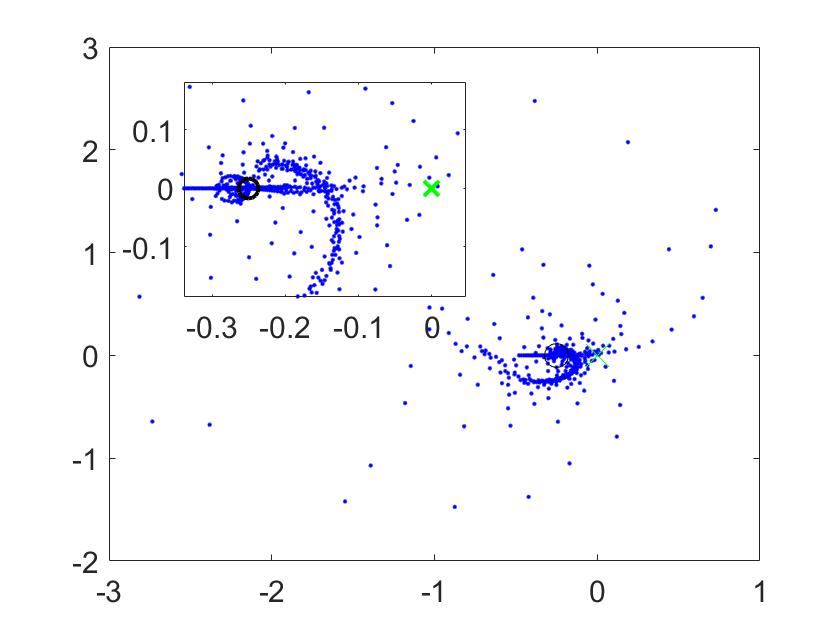} \\
(c) $N^wS^w$, $\lambda=100$ & (d) $\widetilde{N}^wS^w$, $\lambda=100$\\
\includegraphics[scale=0.2]{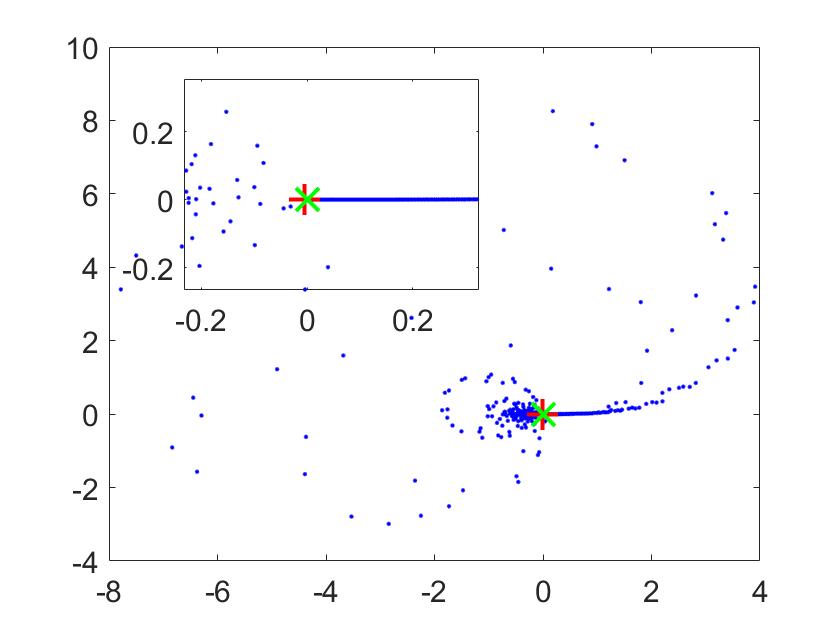} &
\includegraphics[scale=0.2]{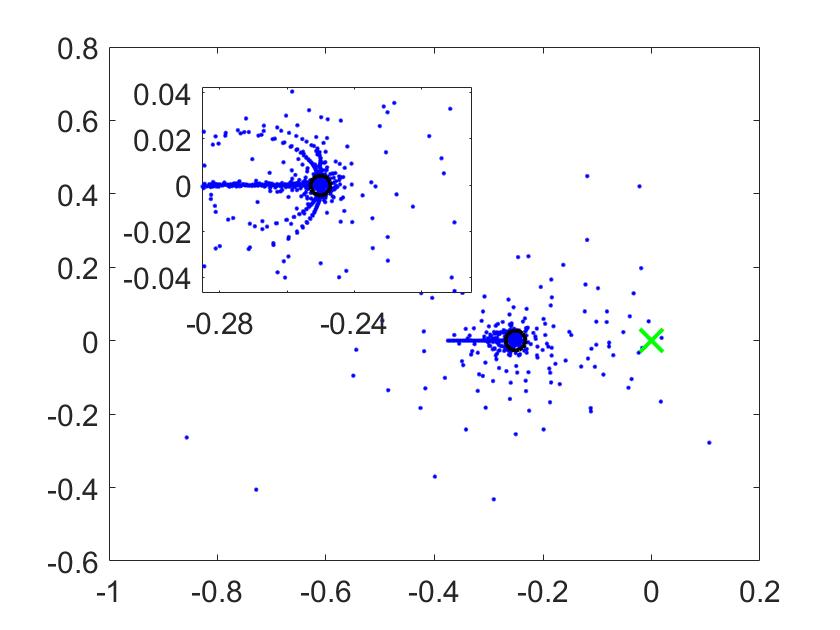} \\
(e) $N^wS^w$, $\lambda=-0.99$ & (f) $\widetilde{N}^wS^w$, $\lambda=-0.99$
\end{tabular}
\caption{Eigenvalue distribution of $N^wS^w$ and $\widetilde{N}^wS^w$ for the spiral-shaped arc. Red points (+): $(-1/4+C_{\lambda,\mu}^2, 0)$; Black
  points ($\circ$): $(-1/4,0)$; Green points ($\times$): $(0,0)$. The
  eigenvalues of $N^wS^w$ and $\widetilde{N}^wS^w$ are bounded away from zero and infinity.}
\label{OpenEeigNS}
\end{figure}

\subsection{Numerical study of the Calder\'{o}n relation for open arcs}
\label{sec:3.3}

Given the smoothness of the solutions of the equations arising from
the weighted operators $S^w$, $N^w$ and $\widetilde{N}^w$, and in
light of equation (\ref{calderonclosed}), it is reasonable to consider
the composite operator $N^wS^w$ as a possible basis for solution of
open-arc problems. The operator $\widetilde{N}^wS^w$ could also
employed in the Dirichlet case. Figure~\eqref{OpenEeigNS} presents the
eigenvalue distribution for discrete versions of the operators
$N^wS^w$ and $\widetilde{N}^wS^w$---which were obtained numerically
for a spiral-shaped open arc using the spectrally accurate quadrature
rules described in section \ref{sec:4.1}. \textcolor{rot}{It can be
  see from Fig.~\ref{OpenEeigNS} that the eigenvalues of $N^wS^w$ and
  $\widetilde{N}^wS^w$ are bounded away from zero and infinity.} A
theoretical analysis of the Calder\'{o}n relation for acoustic
diffraction by open arcs is presented in~\cite{BL12,LB15}. The
corresponding analysis for elastic problems is left for future work.

In view of this numerical evidence and theoretical background on
related problems we suggest that the open arc problems can be solved
effectively by means of the integral equations \be
\label{BIE2}
N^wS^w[\alpha]=N^w[F],\quad
\widetilde{N}^wS^w[\alpha]=\widetilde{N}^w[F],\quad N^wS^w[\beta]=G
\quad\mbox{on}\quad\Gamma.  \en The smooth solutions $\alpha$ and
$\beta$ of these equations are related to the singular solutions of
(\ref{BIEe}) via the relations $\phi=\alpha/w$ and $\psi=\beta\cdot
w$.

\subsection{Regularized formulation for hyper-singular operator}
\label{sec:3.4}

\textcolor{rot}{In this section, we derive an accurate regularized
  formulation for the hyper-singular boundary integral operator
  $\widetilde{N}^w$ in light of the techniques presented
  in~\cite{BXY17,YHX17} for the closed-surface case. The regularized
  formulation for the hyper-singular boundary integral operator $N^w$
  can be obtained directly by letting $\widetilde{\mu}=\mu$.}

\subsubsection{Closed-surface case}

The artificial traction operator (\ref{astress}) can be expressed in the form \be
\label{Tform2}
\widetilde{T}(\pa,\nu)u(x)= (\lambda+\mu)\nu(\nabla \cdot u) + \mu\pa_{\nu}u + \widetilde{\mu} M(\pa,\nu)u
\en
where the operator $M(\pa,\nu)$, whose elements are also called the G\"unter derivatives, is defined as
\ben
M(\pa,\nu)u(x)= \pa_{\nu}u -\nu(\nabla\cdot u)+\nu\times \,\curl\,u.
\enn
In fact, the operator $M$ in 2D is the multiplication of the tangential derivative and a constant matrix, i.e.,
\ben
M(\pa,\nu)u(x)=A\frac{du(x)}{ds_x},\quad A=\begin{pmatrix}
0 & -1 \\
1 & 0
\end{pmatrix}.
\enn
Then following the technique in \cite{YHX17} which is part of the conversion of traction and normal derivatives into tangential derivatives, it can be established that
\be
\label{Nce1}
\widetilde{N}[u](x)&=& -\int_\Gamma\left[\rho\omega^2 (\nu_x\nu_y^\top-\nu_x^\top\nu_yI)-\widetilde{\mu}k_s^2\gamma_{k_s}(x,y)J_{\nu_x,\nu_y}- \rho\omega^2\gamma_{k_p}(x,y)\nu_x\nu_y^\top\right]u(y)ds_y\nonumber\\
&+& (\mu+\widetilde{\mu})^2\frac{d}{ds_x}\int_\Gamma AE(x,y)A\frac{du(y)}{ds_y} ds_y +2(\mu+\widetilde{\mu})\frac{d}{ds_x}\int_\Gamma  \gamma
_{k_s}(x,y)\frac{du(y)}{ds_y}ds_y \nonumber\\
&-& (\mu+\widetilde{\mu}) \int_\Gamma  \nu_{x} \nabla _x^\top [\gamma_{k_s}(x,y)-\gamma_{k_p}(x,y)]A\frac{du(y)}{ds_y}ds_y\nonumber\\
&-& (\mu+\widetilde{\mu})\frac{d}{ds_x}\int_\Gamma A\nabla_y [\gamma_{k_s}(x,y)-\gamma_{k_p}(x,y)]\nu_y^\top u(y)ds_y,
\en
where
\ben
J_{\nu_x,\nu_y}=\nu_y\nu_x^\top-\nu_x\nu_y^\top.
\enn
We omit the proof here.

\subsubsection{Open-arc  case}

\textcolor{rot}{For the open-arc case, due to the smooth boundary-vanishing weight $w$ in $\widetilde{N}^w$, the integration-by-parts formula
\ben
\int_\G \frac{dF(y)}{ds_y}w(y)u(y)ds_y=-\int_\G F(y)\frac{d(w(y)u(y))}{ds_y}ds_y
\enn
holds. We thus, can obtain from the steps in \cite{YHX17} the regularized formulation for the hyper-singular boundary integral operator $\widetilde{N}^w$ as follows
\be
\label{Nce2}
\widetilde{N}^w[u](x)&=& -\int_\Gamma\left[\rho\omega^2 (\nu_x\nu_y^\top-\nu_x^\top\nu_yI)-\widetilde{\mu}k_s^2\gamma_{k_s}(x,y)J_{\nu_x,\nu_y}- \rho\omega^2\gamma_{k_p}(x,y)\nu_x\nu_y^\top\right]w(y)u(y)ds_y\nonumber\\
&+& (\mu+\widetilde{\mu})^2\frac{d}{ds_x}\int_\Gamma AE(x,y)A\frac{d(w(y)u(y))}{ds_y} ds_y +2(\mu+\widetilde{\mu})\frac{d}{ds_x}\int_\Gamma  \gamma
_{k_s}(x,y)\frac{d(w(y)u(y))}{ds_y}ds_y \nonumber\\
&-& (\mu+\widetilde{\mu}) \int_\Gamma  \nu_{x} \nabla _x^\top [\gamma_{k_s}(x,y)-\gamma_{k_p}(x,y)]A\frac{d(w(y)u(y))}{ds_y}ds_y\nonumber\\
&-& (\mu+\widetilde{\mu})\frac{d}{ds_x}\int_\Gamma A\nabla_y [\gamma_{k_s}(x,y)-\gamma_{k_p}(x,y)]\nu_y^\top w(y)u(y)ds_y.
\en}

\section{Numerical experiments}
\label{sec:4}
\subsection{Numerical implementations}
\label{sec:4.1}

In this section, we present spectral quadrature rules for the operators $S^w$, $N^w$ and $\widetilde{N}^w$ which give rise to an efficient and accurate solver for the general elastic and thermoelastic open arc diffraction problems.

As indicated in Section~\ref{sec:3.1}, without loss of generality we may use
a smooth parameterization $x(t)=(x_1(t),x_2(t)), t\in[-1,1]$ of
$\Gamma$ satisfying $|x'(t)|=|dx(t)/dt|>0$, we choose the weight $w$ as
$w(x(t))=\sqrt{1-t^2}$. Hence, the operator $S^w$ give rise to \ben
S^w[\phi](t)=\int_{-1}^1
E(x(t),x(\tau))\frac{\phi(x(\tau))}{\sqrt{1-t^2}}|x'(\tau)|d\tau \enn
Introducing the change of variables $t=\cos\theta$ and
$\tau=\cos\theta'$ and defining $\nu_\theta=\nu_{x(\cos\theta)}$, we
obtain the periodic weighted single-layer operator (using the same
notation) \ben S^w[\phi](\theta)=\int_0^\pi
E(x(\cos\theta),x(\cos\theta'))\widetilde{\phi}(\theta')|x'(\cos\theta')|d\theta',
\quad \widetilde{\phi}(\theta')=\phi(x(\cos\theta')).  \enn It follows
from the definition of the fundamental solution and the series
expansions of Bessel functions~\cite{AS72,CK98}, we can obtain the
decomposition \ben
E(x(t),x(\tau))=E_1(t,\tau)\log|t-\tau|+E_2(t,\tau).  \enn where \ben
E_1(t,\tau) &=& -\frac{1}{2\pi\mu}J_0(k_sr(t,\tau))I+ \frac{1}{2\pi\rho\omega^2} \frac{k_sJ_1(k_sr(t,\tau))-k_pJ_1(k_pr(t,\tau))}{r(t,\tau)}I\\
&-&
\frac{1}{2\pi\rho\omega^2}\frac{(x(t)-x(\tau))^\top(x(t)-x(\tau))}{r(t,\tau)^2}
[k_s^2J_2(k_sr(t,\tau))-k_p^2J_2(k_pr(t,\tau))] \enn with
$r(t,\tau)=|x(t)-x(\tau)|$ and \ben
E_2(t,\tau)=E(x(t),x(\tau))-E_1(t,\tau)\log|t-\tau|.  \enn When
$t=\tau$, we have \ben E_1(t,t)=-\frac{1}{2\pi\mu}I+
\frac{k_s^2-k_p^2}{4\pi\rho\omega^2}I, \enn and \ben
E_2(t,t) &=& \frac{i}{4\mu}\left[1+\frac{2i}{\pi}\left(\log\frac{k_s|x'(t)|}{2}+ C_e \right) \right] +\frac{k_s^2-k_p^2}{4\pi\rho\omega^2}\frac{x'(t)^\top x'(t)}{|x'(t)|^2}\\
&-& \frac{i}{4\rho\omega^2}\left\{\frac{k_s^2-k_p^2}{2}
  \left[1+\frac{2i}{\pi}\left(\log|x'(t)|+ C_e \right)-\frac{i}{\pi}
  \right] +
  \frac{i}{\pi}\left(k_s^2\log\frac{k_s}{2}-k_p^2\log\frac{k_p}{2}
  \right) \right\}.  \enn Use of the Chebyshev points
$\left\{\theta_n=\frac{\pi(2n+1)}{2N}\right\}, n=0,1,\cdots,N-1$ gives
rise to a spectrally convergent cosine representation for smooth,
$\pi$-periodic and even function $\widetilde{\phi}$ as \be
\label{periodicexp}
\widetilde{\phi}(\theta)=\sum_{n=0}^{N-1}a_n\cos(n\theta), \quad a_n=\frac{2-\delta_{0n}}{N}\sum_{j=0}^{N-1}\widetilde{\phi}(\theta_j)\cos(n\theta_j).
\en
It is known from the diagonal property of Symm's operator in the cosine basis $e_n(\theta)=\cos(n\theta)$~\cite{BH07,YS88} that
\be
\label{symmpro}
-\frac{1}{2\pi}\int_0^\pi \log|\cos\theta-\cos\theta'|e_n(\theta')d\theta' =\lambda_ne_n(\theta), \quad \lambda_n=\begin{cases}
\frac{\log2}{2}, & n=0, \\
\frac{1}{2n}, & n\ge 1.
\end{cases}
\en
Applying equation \ref{symmpro} to each term of expansion (\ref{periodicexp}) we can obtain the well-known spectral quadrature rule for the logarithmic kernel
\be
\label{sqr}
\int_0^\pi \log|\cos\theta-\cos\theta'|\widetilde{\phi}(\theta')d\theta' \sim \frac{\pi}{N}\sum_{j=0}^{N-1}\widetilde{\phi}(\theta_j)R_j^{(N)}(\theta),
\en
where
\be
\label{RjN}
R_j^{(N)}(\theta)=-2\sum_{m=0}^{N-1}(2-\delta_{0m}) \lambda_m\cos(m\theta_j)\cos(m\theta).
\en
Together with the trapezoidal integration for smooth function we therefore obtain the spectrally quadrature approximation of the operator $S_{o}$ that
\be
\label{sqaS}
S^w[\phi](\theta)\sim \frac{\pi}{N}\sum_{j=0}^{N-1} \widetilde{\phi}(\theta_j)|x'(\cos\theta_j)| \left[E_1(\cos\theta,\cos\theta_j)R_j^{(N)}(\theta) + E_2(\cos\theta,\cos\theta_j)\right].
\en
Then we can evaluate $S^w[\phi]$ in the sets of quadrature points $\{\theta_n, n=0,\cdots,N-1\}$ by means of a matrix-vector multiplication involving the matrix ${\bf S}^w$ whose elements are defined by
\ben
[{\bf S}^w]_{ij}=\frac{\pi}{N}|x'(\cos\theta_j)| \left[E_1(\cos\theta_i,\cos\theta_j)R_j^{(N)}(\theta_i) + E_2(\cos\theta_i,\cos\theta_j)\right],
\enn
in which the quantities $R_j^{(N)}(\theta_i)$ can be evaluated through
\ben
R_j^{(N)}(\theta_i)=R^{(N)}(|i-j|)+R^{(N)}(i+j+1),
\enn
where
\ben
R^{(N)}(l)=-\sum_{m=0}^{N-1} (2-\delta_{0m})\lambda_m\cos\left(\frac{lm\pi}{N}\right), \quad l=0,\cdots,2N-1,
\enn
while this expression can significantly reduce the computational cost of ${\bf S}_{o}^w$.

In order to evaluate $N^w$ and $\widetilde{N}^w$, we use the regularized formulations given in section \ref{sec:3.4}. In (\ref{Nce1}), the first term on the right hand side can be evaluated by means of a rule analogous to the single-layer operator of acoustic~\cite{BL12}). The other terms take one of the following forms
\be
\label{otherterm}
D_0S_1T_0, \quad S_2T_0, \quad D_0S_3,
\en
where $S_j,j=1,2,3$, whose kernels are of logarithmic type, can be evaluated analogous to $S^w$, and
\be
D_0[\widetilde{\phi}](\theta)=\frac{1}{\sin\theta} \frac{d\widetilde{\phi}(\theta)}{d\theta}=\frac{d\phi(t)}{dt},\quad T_0[\widetilde{\phi}](\theta)= \frac{d}{d\theta}(\widetilde{\phi}(\theta)\sin\theta).
\en
Then we approximate the quantity $T_0[\widetilde{\phi}]$ by means of term per term differentiation of the sine expansion of the function $\widetilde{\phi}(\theta)\sin\theta$ (which can itself be produced efficiently by means of an FFT). The quantity $D_0[\widetilde{\phi}]$ can be evaluated by invoking classical FFT-based Chebyshev differentiation rules~\cite{PTVF92}).

\subsection{Numerical examples}
\label{sec:4.2}

The numerical examples presented in what follows were obtained by means of a Matlab implementation of the quadrature rules introduced in section \ref{sec:4.1} for numerical evaluation of the operators $S^w$, $N^w$ and $\widetilde{N}^w$ in conjunction with the iterative linear algebra solver GMRES. In all cases the (maximum) errors reported were evalueated by comparisons with exact or highly-resolved numerical solutions. Select the elastic parameters as $\rho=1$, $\mu=1$ and $\lambda=2$. We call the equations in (\ref{BIE1}) and (\ref{BIE2}) as $\mbox{Dir}(S^w)$, $\mbox{Neu}(N^w)$ and $\mbox{Dir}(N^wS^w)$, $\mbox{Dir}(\widetilde{N}^wS^w)$, $\mbox{Neu}(N^wS^w)$, respectively.

\begin{table}[ht]
\caption{Near-field errors for elastic scattering by a unit circle. GMRES tol: $10^{-12}$.}
\centering
\begin{tabular}{ccccccc}
\hline
$\omega$ & $N$ & $\mbox{Dir}(S)$ & $\mbox{Dir}(NS)$ & $\mbox{Dir}(\widetilde{N}S)$ & $\mbox{Neu}(N)$ & $\mbox{Neu}(NS)$ \\
\hline
    & 30  & $6.59\times10^{-4}$  & $6.59\times10^{-4}$  & $6.59\times10^{-4}$ & $9.53\times10^{-3}$  & $9.53\times10^{-3}$ \\
 10 & 40  & $1.55\times10^{-6}$  & $1.55\times10^{-6}$  & $1.55\times10^{-6}$ & $6.02\times10^{-5}$  & $6.02\times10^{-5}$ \\
    & 60  & $2.54\times10^{-12}$ & $2.50\times10^{-12}$ & $2.68\times10^{-12}$ & $2.13\times10^{-11}$ & $2.16\times10^{-11}$ \\
\hline
    & 150  & $1.33\times10^{-7}$  & $1.33\times10^{-7}$  & $1.33\times10^{-7}$ & $1.58\times10^{-7}$  & $1.58\times10^{-7}$ \\
 50 & 160  & $2.51\times10^{-10}$  & $2.51\times10^{-10}$  & $2.54\times10^{-10}$ & $1.27\times10^{-9}$  & $1.27\times10^{-9}$ \\
    & 200  & $7.12\times10^{-13}$ & $1.58\times10^{-12}$ & $9.53\times10^{-12}$ & $2.22\times10^{-12}$ & $1.65\times10^{-12}$ \\
\hline
\end{tabular}
\label{Table1}
\end{table}

We first demonstrate the accuracy of the regularized formulation given in section \ref{sec:3.4}. To see this, let $\Gamma$ be a unit circle and the exact solutions for elastic and thermoelastic problems are given by $u=\nabla_xH_0^{(1)}(k_p|x-z_0|)$ and $u=E_{12}(x,z_0)$, $p=E_{22}(x,z_0)$, respectively where $z_0=(0,0.5)$. We apply the Nystr\"om method \cite{CK98} for the discretization of the operators $S$, $N$ and $\widetilde{N}$. Table 1 shows the spectral (exponentially-fast) convergence which further demonstrate the accuracy of our regularized formulations.

\begin{table}[ht]
\caption{Near-field errors for elastic scattering by a Spiral-Shaped Arc. GMRES tol: $10^{-8}$.}
\centering
\begin{tabular}{ccccccc}
\hline
$\omega$ & $N$ & $\mbox{Dir}(S^w)$ & $\mbox{Dir}(N^wS^w)$ & $\mbox{Dir}(\widetilde{N}^wS^w)$ & $\mbox{Neu}(N^w)$ & $\mbox{Neu}(N^wS^w)$ \\
\hline
    & 100  & $2.97\times10^{-2}$  & $2.97\times10^{-2}$  & $2.97\times10^{-2}$ & $2.22\times10^{-1}$  & $2.22\times10^{-1}$ \\
 10 & 150  & $7.41\times10^{-5}$  & $7.41\times10^{-5}$  & $7.41\times10^{-5}$ & $4.33\times10^{-3}$  & $4.33\times10^{-3}$ \\
    & 200  & $2.06\times10^{-11}$ & $2.39\times10^{-11}$ & $2.10\times10^{-11}$ & $7.42\times10^{-10}$ & $7.43\times10^{-10}$ \\
\hline
    & 600  & $6.11\times10^{-3}$  & $1.95\times10^{-4}$  & $6.11\times10^{-3}$ & $1.49\times10^{-1}$  & $4.22\times10^{-2}$ \\
 50 & 800  & $4.76\times10^{-10}$  & $8.51\times10^{-8}$  & $1.19\times10^{-7}$ & $1.46\times10^{-6}$  & $1.86\times10^{-6}$ \\
    & 1000 & $2.01\times10^{-14}$ & $1.33\times10^{-7}$ & $1.97\times10^{-12}$ & $1.57\times10^{-8}$  & $9.36\times10^{-9}$ \\
\hline
\end{tabular}
\label{Table2}
\end{table}

Now we consider the elastic scattering by open arcs. Let $u^{inc}$ be a plane incident pressure wave $u^{inc}=d_{inc}\exp(ik_px\cdot d_{inc})$ where $d_{inc}=(\cos\theta_{inc},\sin\theta_{inc})$ and $\theta_{inc}$ is the incident direction. This means that the boundary data is given by
\ben
F=-u^{inc},\quad G=-T(\pa,\nu)u^{inc}.
\enn
To demonstrated the high-order character of the algorithm we consider the elastic scattering by the 'Spiral-Shaped Arc' characterized by the parameterization
\ben
x(t)=\exp(t)(\cos 5t,\sin 5t),\quad t\in[-1,1].
\enn
Choose $\theta_{inc}=\pi/4$. We present the numerical errors of the solution produced by means of the operators $S^w$, $N^w$, $\widetilde{N}^wS^w$ and $N^wS^w$ for the elastic scattering problems which demonstrate the spectral (exponentially-fast) convergence.

\begin{table}[ht]
\caption{Iterations for elastic scattering by a Flat Strip. GMRES tol: $10^{-5}$.}
\centering
\begin{tabular}{ccccccc}
\hline
$\omega$ & $N$ & $\mbox{Dir}(S^w)$ & $\mbox{Dir}(N^wS^w)$ & $\mbox{Dir}(\widetilde{N}^wS^w)$ & $\mbox{Neu}(N^w)$ & $\mbox{Neu}(N^wS^w)$ \\
\hline
10  & 160  & 21 & 12 & 12 & 31 & 13 \\
30  & 480  & 33 & 23 & 16 & 88 & 24 \\
50  & 800  & 41 & 30 & 19 & 137 & 34 \\
80  & 1280 & 48 & 43 & 21 & 204 & 47 \\
100 & 1600 & 51 & 50 & 23 & 246& 58 \\
\hline

\end{tabular}
\label{Table3}
\end{table}

\begin{table}[ht]
\caption{Computing time (seconds) required for elastic scattering by a Flat Strip. GMRES tol: $10^{-5}$.}
\centering
\begin{tabular}{ccccccc}
\hline
$\omega$ & $N$ & $\mbox{Dir}(S^w)$ & $\mbox{Dir}(N^wS^w)$ & $\mbox{Dir}(\widetilde{N}^wS^w)$ & $\mbox{Neu}(N^w)$ & $\mbox{Neu}(N^wS^w)$ \\
\hline
10  & 160  & 0.01 & 0.06 & 0.06 & 0.17 & 0.06 \\
30  & 480  & 0.07 & 0.41 & 0.29 & 2.9 & 0.53 \\
50  & 800  & 0.22 & 1.3 & 0.73 & 11.2 & 1.9 \\
80  & 1280 & 0.60 & 4.4 & 2.1 & 40.4 & 8.9 \\
100 & 1600 & 0.97 & 7.7 & 3.3 & 76.3 & 15.5 \\
\hline

\end{tabular}
\label{Table4}
\end{table}

\begin{table}[ht]
\caption{Iterations for elastic scattering by a Spiral-Shaped Arc. GMRES tol: $10^{-5}$.}
\centering
\begin{tabular}{ccccccc}
\hline
$\omega$ & $N$ & $\mbox{Dir}(S^w)$ & $\mbox{Dir}(N^wS^w)$ & $\mbox{Dir}(\widetilde{N}^wS^w)$ & $\mbox{Neu}(N^w)$ & $\mbox{Neu}(N^wS^w)$ \\
\hline
10  & 160  & 73  & 75  & 61  & 198 & 76 \\
30  & 480  & 109 & 139 & 94  & 444 & 145 \\
50  & 800  & 126 & 196 & 113 & 651 & 204 \\
80  & 1280 & 147 & 274 & 139 & 985 & 291 \\
100 & 1600 & 160 & 330 & 157 & 1222& 348 \\
\hline

\end{tabular}
\label{Table5}
\end{table}

\begin{table}[ht]
\caption{Computing time (seconds) for elastic scattering by a Spiral-Shaped Arc. GMRES tol: $10^{-5}$.}
\centering
\begin{tabular}{ccccccc}
\hline
$\omega$ & $N$ & $\mbox{Dir}(S^w)$ & $\mbox{Dir}(N^wS^w)$ & $\mbox{Dir}(\widetilde{N}^wS^w)$ & $\mbox{Neu}(N^w)$ & $\mbox{Neu}(N^wS^w)$ \\
\hline
10  & 160  & 0.08 & 0.25 & 0.18 & 0.8   & 0.23 \\
30  & 480  & 0.36 & 2.7 & 1.6 & 10.7  & 3.7 \\
50  & 800  & 0.93 & 15.9 & 4.6 & 41.9  & 15.7 \\
80  & 1280 & 2.5  & 47.6 & 13.7 & 183.2 & 51.2 \\
100 & 1600 & 3.7 & 53.9 & 23.7 & 347.8 & 93.8 \\
\hline

\end{tabular}
\label{Table6}
\end{table}

In Table 3 and 4, we present the iteration numbers required to achieve the GMRES tolerance $10^{-5}$ for the scattering by a flat strip $[-1,1]$ and the exponential spiral mentioned above. It can be seen that for Neumann problems, $\mbox{Neu}(N^w)$ requires very large number of iterations as the frequency grows and, thus, the computing times required by the low-iteration equation $\mbox{Neu}(N^wS^w)$ are significantly lower than those required by $\mbox{Neu}(N^w)$. But for the Dirichlet problem, $\mbox{Dir}(N^wS^w)$ requires almost the same number of iterations as $\mbox{Neu}(N^wS^w)$ which is significantly larger that $\mbox{Dir}(S^w)$ for the spiral-shaped arc. But the corresponding equation $\mbox{Dir}(\widetilde{N}^wS^w)$ requires fewer iterations than $\mbox{Dir}(S^w)$ for both flat strip and spiral-shaped arc. We present the corresponding total time required by the solver to find the solution once the needed matrixes for iteration, $S^w$, $S_1$, $S_2$ and $S_3$ mentioned in section ?, are stored. It should be pointed out that the total computing cost of the equations $\mbox{Dir}(\widetilde{N}^wS^w)$ and $\mbox{Dir}(N^wS^w)$ is generally higher than $\mbox{Dir}(S^w)$ since the application of the operator in $\mbox{Dir}(S^w)$ is significantly less expensive than the application of operator in $\mbox{Dir}(\widetilde{N}^wS^w)$ and $\mbox{Dir}(N^wS^w)$. In addition, from the regularized formulation, the cost for the generation of the matrixes $\widetilde{N}^w$ and $N^w$ (except the computing of tangential derivatives) is about 5 times of that for $S^w$.

\begin{figure}[htbp]
\centering
\begin{tabular}{ccc}
\includegraphics[scale=0.25]{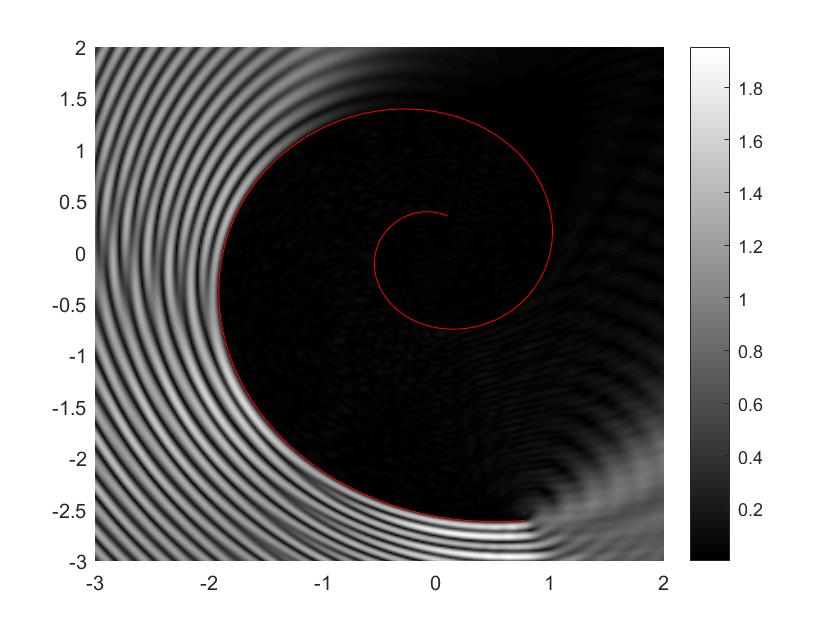} &
\includegraphics[scale=0.25]{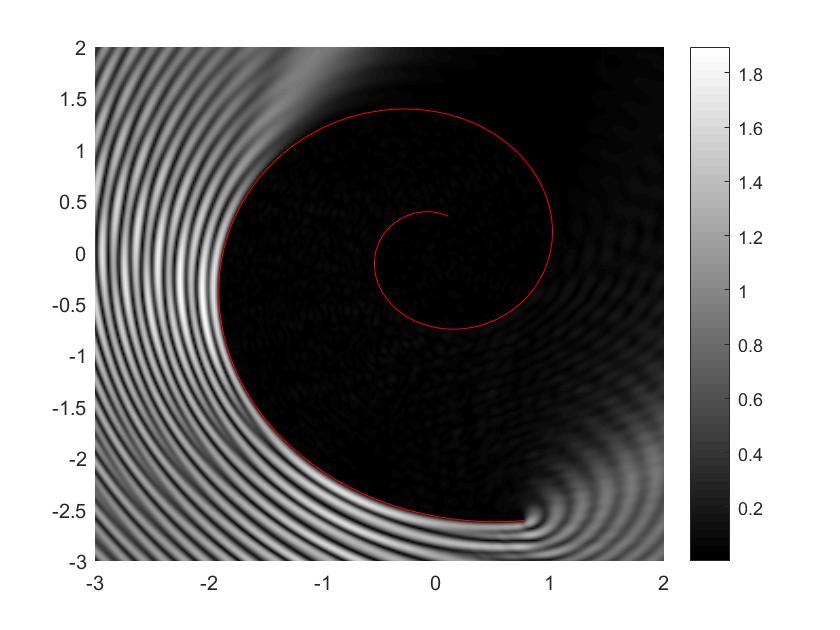} \\
(a) $|u_1|, \theta_{inc}=\pi/4$ & (b) $|u_2|, \theta_{inc}=\pi/4$  \\
\includegraphics[scale=0.25]{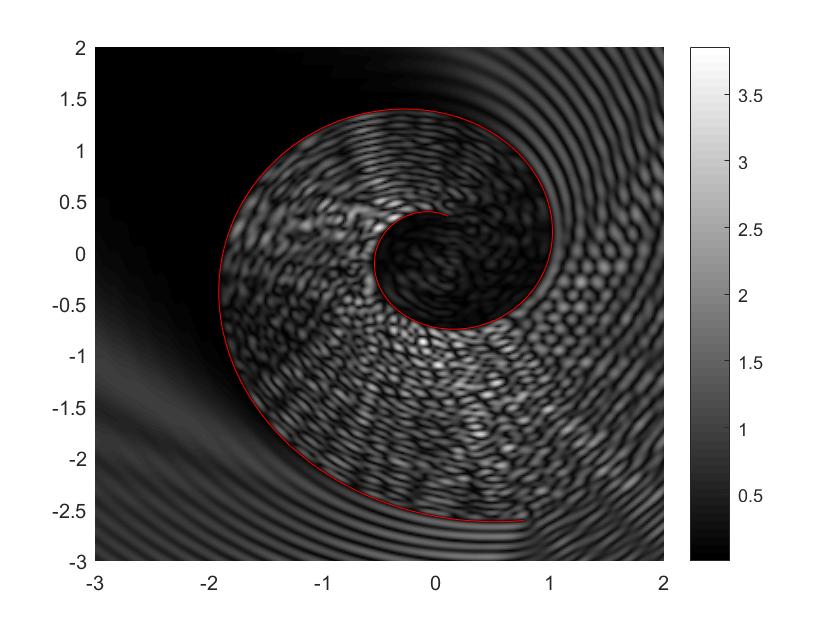} &
\includegraphics[scale=0.25]{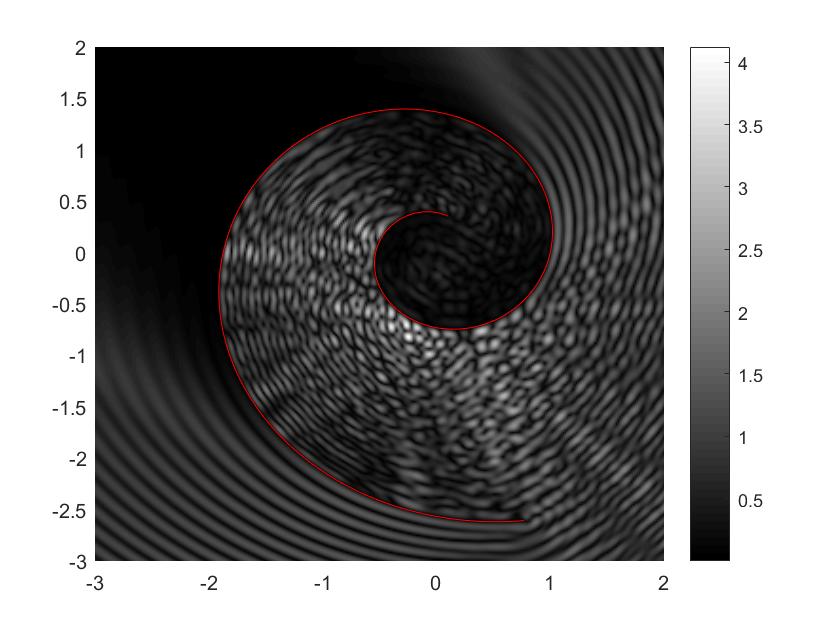} \\
(c) $|u_1|, \theta_{inc}=3\pi/4$ & (d) $|u_2|, \theta_{inc}=3\pi/4$
\end{tabular}
\caption{Elastic scattering by a Spiral-Shaped Arc with Dirichlet boundary condition. GMRES tol: $10^{-5}$.}
\label{FigESol12}
\end{figure}

\begin{figure}[htbp]
\centering
\begin{tabular}{ccc}
\includegraphics[scale=0.25]{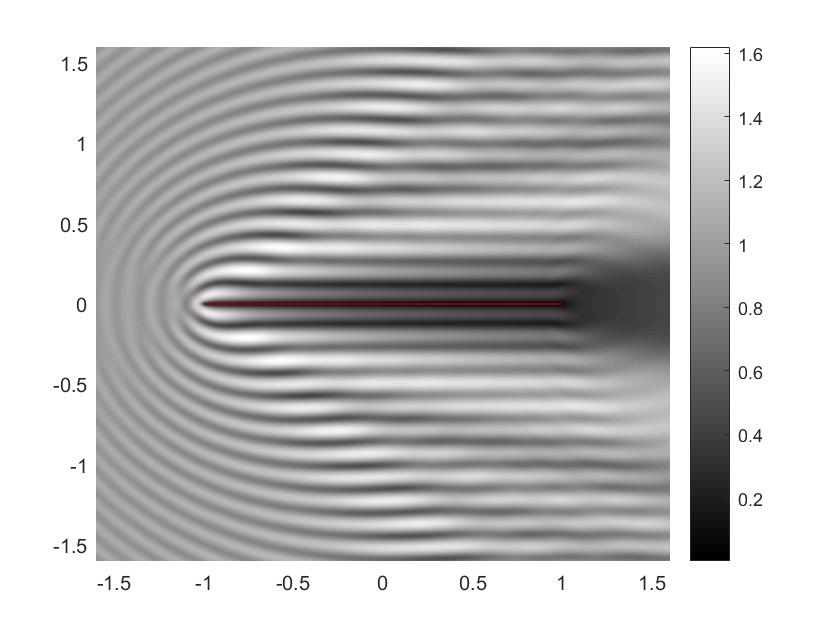} &
\includegraphics[scale=0.25]{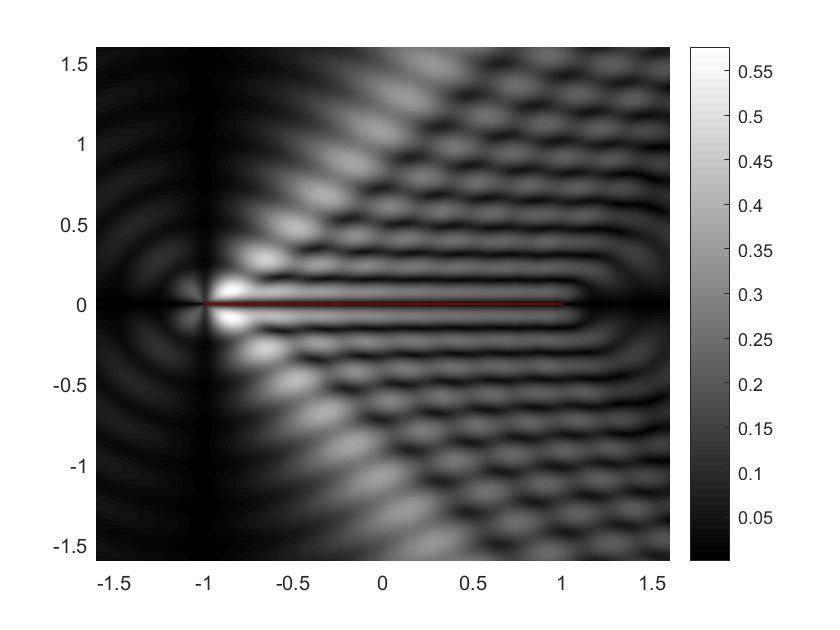} \\
(a) $|u_1|, \theta_{inc}=0$ & (b) $|u_2|, \theta_{inc}=0$  \\
\includegraphics[scale=0.25]{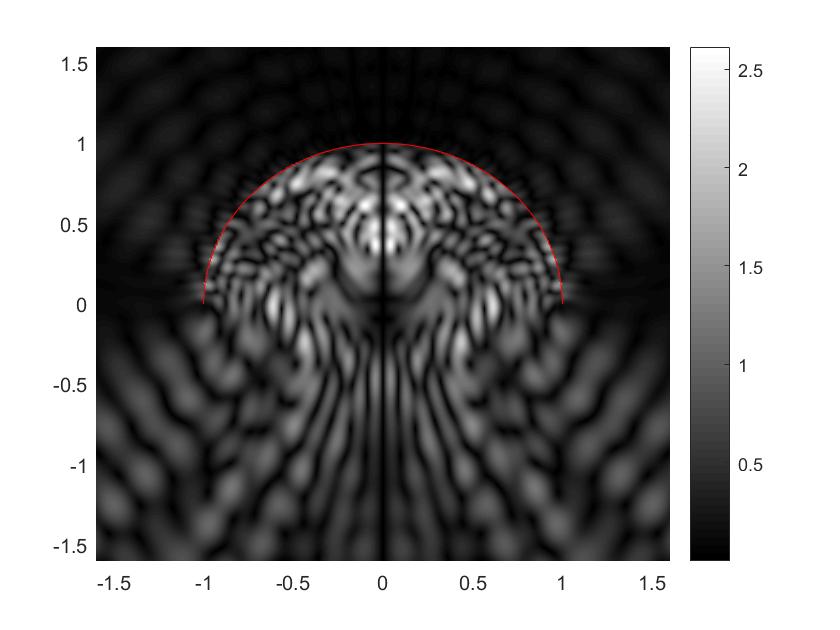} &
\includegraphics[scale=0.25]{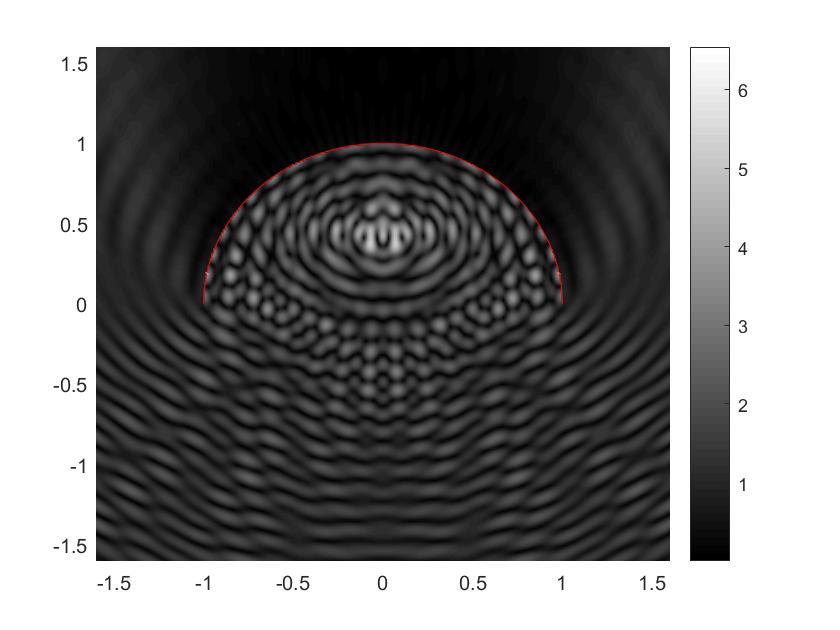} \\
(c) $|u_1|, \theta_{inc}=\pi/2$ & (d) $|u_2|, \theta_{inc}=\pi/2$
\end{tabular}
\caption{Elastic diffraction patterns with Dirichlet condition (a,b) and Neumann condition (c,d). GMRES tol: $10^{-5}$.}
\label{FigESol34}
\end{figure}

In order to obtain an indication of the manner in which an open arc problem can be viewed as a limit of closed-curve problems and provide an independent verification of the validity of our solvers, we consider a test case in which the flat strip $[-1,1]$ is viewed as the limit as $a\rightarrow0$ of the family of closed curves $x(t)=(\cos t,a\sin t)$. It is known that the scattered field $u$ admits the asymptotic behavior
\ben
u(x)=\frac{e^{ik_px+i\pi/4}}{\sqrt{8\pi k_p|x|}}u_p^\infty(\hat{x})\hat{x} +\frac{e^{ik_sx+i\pi/4}}{\sqrt{8\pi k_s|x|}}u_s^\infty(\hat{x})\hat{x}^\perp +\mathcal{O}(|x|^{-3/2}),\quad |x|\rightarrow\infty.
\enn
Then for the Dirichlet problem, the P-part $u_p^\infty$ and S-part $u_p^\infty$ of the far-field pattern of $u$ are given by
\ben
u_p^\infty &=& \int_\Gamma e^{-ik_p\hat{x}\cdot y}[\hat{x}\cdot\phi(y)]/w(y)ds_y, \\
u_s^\infty &=& \int_\Gamma e^{-ik_s\hat{x}\cdot y}[\hat{x}^\perp\cdot\phi(y)]/w(y)ds_y
\enn
Using the closed-curve Nystr\"om method we evaluate the scattered field of Dirichlet problems for values of $a$ approaching 0 and present the P-part and S-part of the far-field pattern side-by-side the corresponding far-field pattern for the limiting open arc as produced by the $S^w$-based open-arc solver. Clearly, the closed-curve and open-arc solutions are quite close to each other.

\begin{figure}[htbp]
\centering
\begin{tabular}{cc}
\includegraphics[scale=0.25]{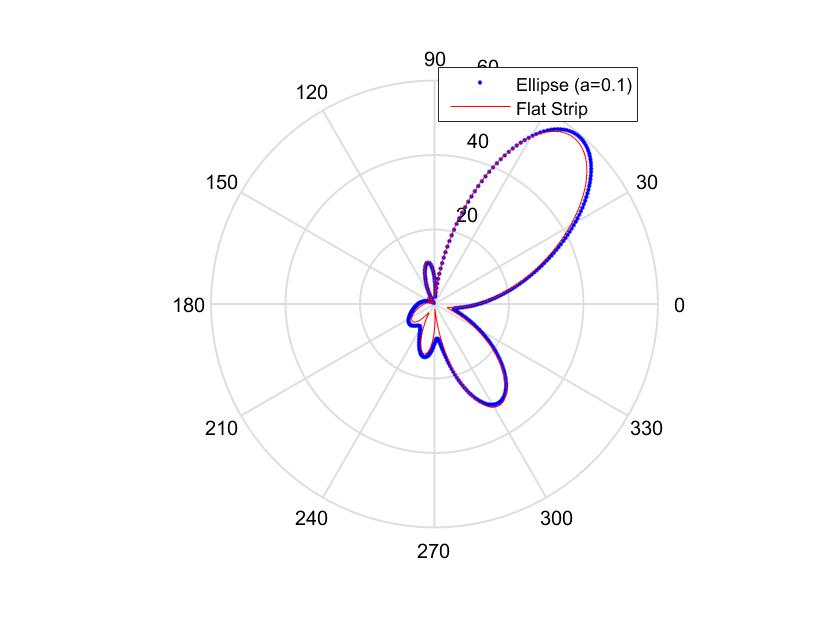} &
\includegraphics[scale=0.25]{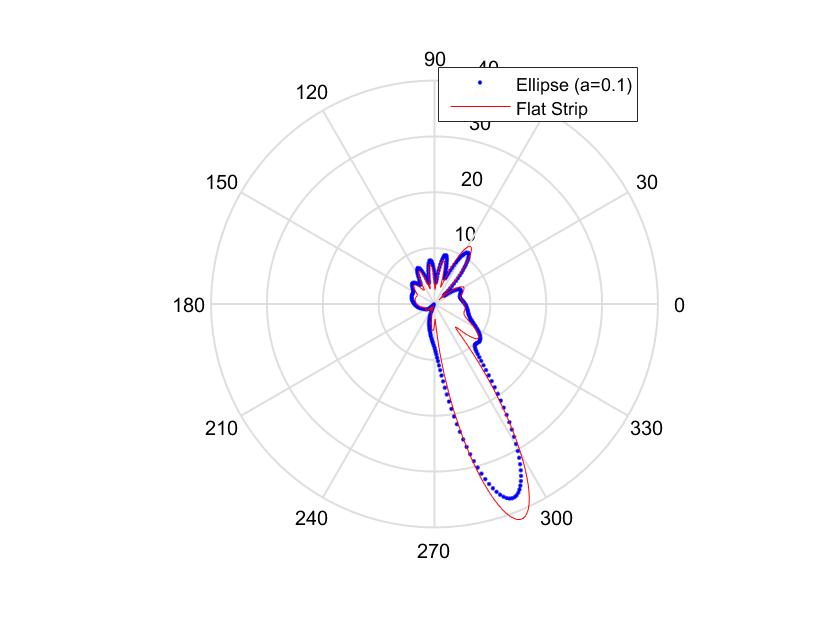} \\
(a) & (b) \\
\includegraphics[scale=0.25]{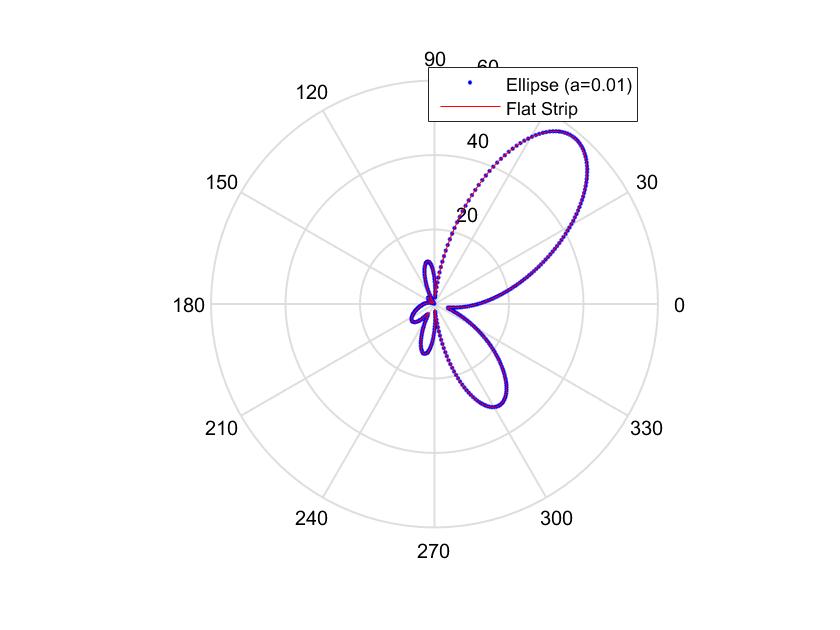} &
\includegraphics[scale=0.25]{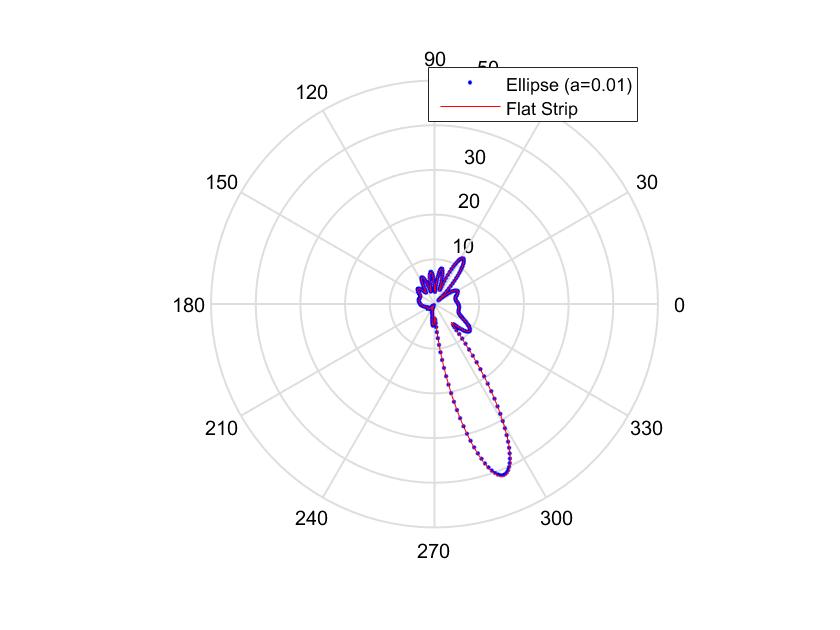} \\
(c) & (d) \\
\includegraphics[scale=0.25]{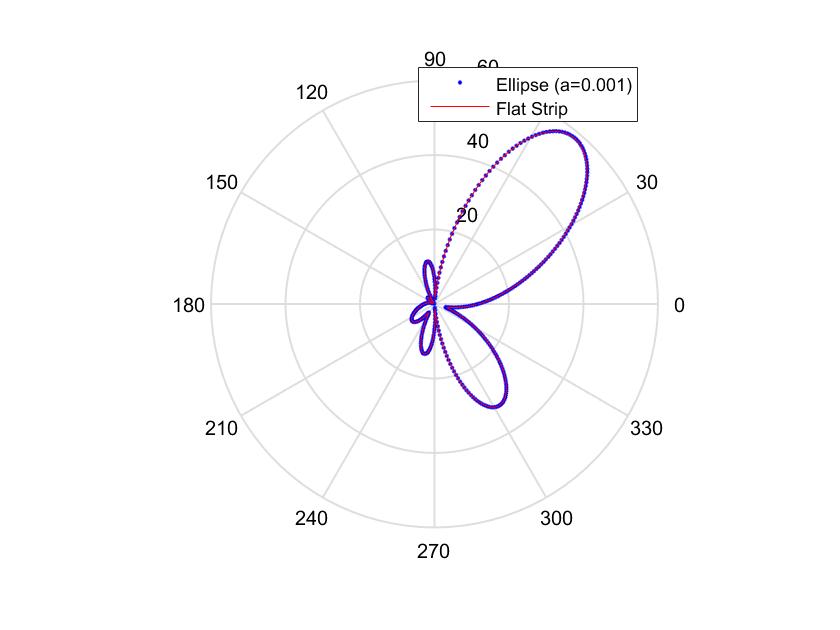} &
\includegraphics[scale=0.25]{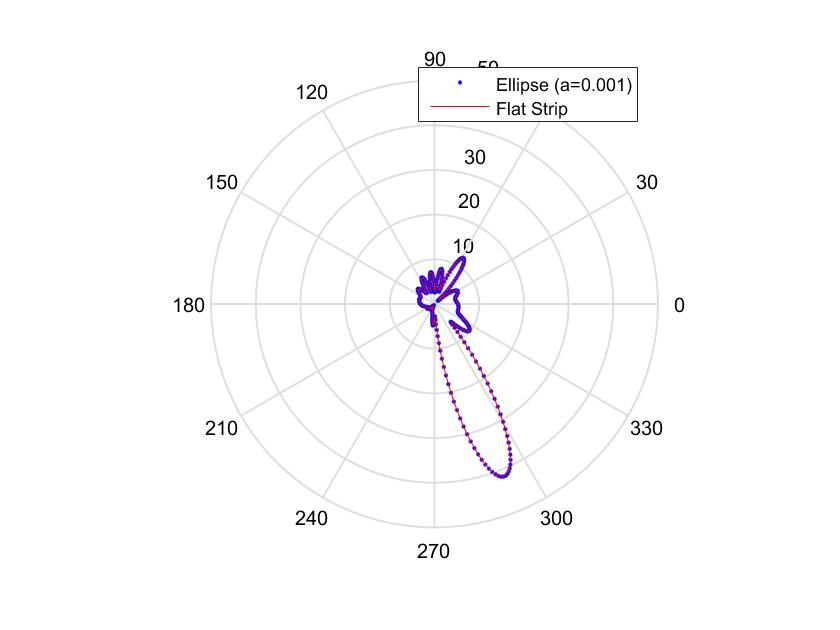} \\
(e) & (f) \\
\end{tabular}
\caption{Far-field patterns $u_p^\infty$ (a,c,e) and $u_s^\infty$ (b,d,f) for the scattering by a sequence of increasingly thin closed ellipse converging to the flat strip.}
\label{FarField}
\end{figure}

\section{Conclusion}
\label{sec:5}

We have introduced new integral solvers and associated numerical algorithms for the elastic scattering by open arcs with Dirichlet or Neumann boundary condition in two dimensions. The new methods enjoy spectral convergence and reduce the number of GMRES iterations consistently across various geometries and frequency regimes. In particular, the new formulation is highly beneficial for the Neumann problem, giving rise to order-of-magnitude improvements in computing time over the original hypersingular formulation. Theoretical investigation of the Calder\'on formula in open-arc case and the generalization of the method enabling efficient solution of problems of elastic and thermoelastic scattering by open-surfaces in three dimensions are left for future work.

\section*{Acknowledgments}

OB gratefully acknowledge support by NSF, AFOSR and DARPA through
contracts DMS- 1411876 and FA9550-15-1-0043 and HR00111720035, and the
NSSEFF Vannevar Bush Fellowship under contract number
N00014-16-1-2808. LX is partially supported by a Key Project of the
Major Research Plan of NSFC (No. 91630205), and a NSFC Grant
(No. 11771068).

\end{document}